\documentclass[10pt]{amsart}
\usepackage{amssymb,latexsym,amsmath,amsfonts}

\numberwithin{equation}{section}
\theoremstyle{definition}
\newtheorem{definition}{Definition}[section]

\theoremstyle{plain}
\newtheorem{theorem}[definition]{Theorem}
\newtheorem{lemma}[definition]{Lemma}

\newcommand{\Cn}{{\mathbb{C}}^n}
\newcommand{\RN}{{\mathbb{R}}^N}

\newcommand{\bdy}{\partial \Omega}

\newcommand{\OM}{\Omega}
\newcommand{\mn}{\boldsymbol{M}}
\newcommand{\strat}{\mathcal{M}}
\newcommand{\st}{\widetilde{S}}
\newcommand{\STrat}{\mathcal{S}}
\newcommand{\varST}{\widetilde{\mathcal{S}}}
\newcommand{\Nu}{\mathfrak{N}}
\newcommand{\er}{\text{Re}}

\newcommand{\vardrv}[2]{\sum_{j=1}^n\partial_j\rho(#1)(#2)}
\newcommand{\drv}[2]{\sum_{j=1}^n\partial_j\rho(#1)[#2]}
\newcommand{\dbdrv}{\partial^2_{jk}\rho}
\newcommand{\dbdrvc}{\partial^2_{j\overline{k}}\rho}
\newcommand{\derrho}[1]{\partial_j\rho(#1)}

\newcommand{\zbar}{\overline{z}}

\newcommand{\Rd}{{\mathbb{R}}^d}
\newcommand{\gamder}[2]{\dfrac{\partial{\gamma}_{#1}}{\partial x_{#2}}}
\newcommand{\fuder}{\dfrac{\partial^2\gamma_j}{\partial x_{\mu}\partial x_{\nu}}(y)}
\newcommand{\zum}{\sum_{\mu=1}^d}
\newcommand{\zuma}{\sum_{\mu,\nu=1}^d}
\newcommand{\add}{\sum_{j=1}^n}
\newcommand{\adda}{\sum_{j,k=1}^n}

\newcommand{\eps}{\varepsilon}

\newcommand{\zt}{\zeta}

\newcommand{\form}{\mathfrak{H}}
\newcommand{\yz}{y^z_r}

\newcommand{\BDO}{B_d}
\newcommand{\BDC}{\mathbb{B}_d}

\newcommand{\Dim}{\text{dim}_{\mathbb{R}}}
\newcommand{\Rnk}{\text{rank}_{\mathbb{R}}}
\newcommand{\tansp}[2]{T_{#2}[\gamma(S_{#1,\alpha})]}
\newcommand{\umm}{\mathcal{M}^\alpha_{d-m-1}}

\newcommand{\rod}{r_*}

\newcommand{\bola}{\boldsymbol{A}\boldsymbol{(}\boldsymbol{\Omega}\boldsymbol{)}}
\newcommand{\Can}{\boldsymbol{\mathbb{C}}^{\boldsymbol{n}}}
\newcommand{\Q}{\sum_{j,k=1}^n\left(\partial^2_{j\overline{k}}\rho(\gamma(x)) \ 
		[d\gamma(x)v]_j \ \overline{[d\gamma(x)v]}_k+
		\partial^2_{jk}\rho(\gamma(x)) \ [d\gamma(x)v]_j \ [d\gamma(x)v]_k\right)}
\newcommand{\varQ}{\sum_{j,k=1}^n\left(\partial^2_{j\overline{k}}\rho(\gamma(s)) \ 
		[d\gamma(s)v]_j \ \overline{[d\gamma(s)v]}_k+
		\partial^2_{jk}\rho(\gamma(s)) \ [d\gamma(s)v]_j \ [d\gamma(s)v]_k\right)}

\begin{document}
\title[Peak-interpolation manifolds in convex domains]{On
Peak-interpolation Manifolds for \\
$\bola$ for Convex Domains in $\Can$}
\author{Gautam Bharali}
\address{Mathematics Department, University of Wisconsin-Madison, 480 Lincoln Drive,
Madison, WI 53706}
\email{bharali@math.wisc.edu}

\begin{abstract} Let $\Omega$ be a bounded, weakly convex domain in
$\Cn$, $n\geq 2$, having real-analytic boundary. $A(\OM)$ is the algebra of all functions 
holomorphic in $\OM$ and continuous upto the boundary. A submanifold $\mn\subset\bdy$ is said
to be complex-tangential if $T_p(\mn)$ lies in the maximal complex subspace of $T_p(\bdy)$ for
each $p\in\mn$. We show that for real-analytic submanifolds $\mn\subset\bdy$, if $\mn$ is 
complex-tangential, then every compact subset of $\mn$ is a peak-interpolation set for $A(\OM)$.
\end{abstract}

\maketitle

\section{Statement of Main Result}

Let $\OM $ be a bounded domain in $\Cn $, and let $A(\OM )$ be the algebra of
functions holomorphic in $\OM $ and continuous upto the boundary. Recall that a compact
subset $K \subset \bdy $ is called a {\bf peak-interpolation set for $\bola$} if given
any $f \in \mathcal{C}(K)$, $f\not\equiv 0$, there exists a function $F \in A(\OM )$ such that
$F|_K = f$ and $|F(\zeta)| < \sup_K|f|$ for every $\zeta \in \overline{\OM}\setminus K$.
\medskip

We are interested in determining when a sufficiently smooth submanifold $\mn\subset\bdy$ is a 
peak-interpolation set for $A(\OM)$. When $\OM$ is a strictly pseudoconvex domain having $\mathcal{C}^2$ 
boundary, and $\mn$ is of class $\mathcal{C}^2$, the situation is very well understood; refer to the works 
of Henkin \& Tumanov \cite{HT}, Nagel \cite{aN}, and Rudin \cite{wR}. In the strictly pseudoconvex setting, 
$\mn$ is a peak-interpolation set for $A(\OM)$ if and only if $\mn$ is {\bf complex-tangential}, i.e.
$T_p(\mn)\subset H_p(\bdy) \ \forall p\in\mn$. Here, and in what follows, for any submanifold 
$\mathcal{M}\subseteq\bdy$, $T_p(\mathcal{M})$ will denote the real tangent space to $\mathcal{M}$ at the
point $p\in\mathcal{M}$, while $H_p(\bdy)$ will denote the maximal complex subspace of $T_p(\bdy)$.
\medskip

Very little is known, however, when $\OM$ is a weakly pseudoconvex of finite type (There are several notions
of type for domains in $\Cn, \ n\geq 3$. We shall not define them at this juncture; the interested reader
may refer to \cite{BG}, \cite{dC}, \cite{jdA}, \cite{McN}.). In view of a result by Nagel \& Rudin
\cite{NR}, it is still necessary for $\mn$ to be complex-tangential. However, showing even that any smooth,
(topologically) closed complex-tangential arc in $\bdy$ is a peak-interpolation set for $A(\OM)$, for a
general smoothly bounded weakly pseudoconvex domain of finite type, is a difficult problem. This is 
because doing so would necessarily imply that every point in $\bdy$ is a peak point for $A(\OM)$. Whether 
or not this is true for general pseudoconvex domains of finite type is an extremely difficult open question in 
the theory of functions in several complex variables. In this paper we show that when $\OM$ is a 
{\em convex} domain and $\bdy$ and $\mn$ are real-analytic, it suffices for $\mn$ to be complex-tangential for
it to be a peak-interpolation set for $A(\OM)$.
\medskip

\pagebreak
 
Our main result is as follows :

\begin{theorem}\label{T:main} Let $\OM$ be a bounded (weakly) convex domain in $\Cn, \ n\geq 2$, having
real-analytic boundary, and let $\mn$ be a real-analytic submanifold of $\bdy$. If $\mn$ is complex-tangential,
then $\mn$ (and thus, every compact subset of $\mn$) is a peak-interpolation set for $A(\OM)$.
\end{theorem}
\medskip

\section{Some notation and introductory remarks}

In what follows, the notation $\langle \ , \ \rangle$ will denote the usual real inner product
on $\Rd$. Furthermore, given vectors $v,w\in\Rd$ and a real $d\times d$ matrix $M=[m_{jk}]$, the notation
$\langle v \ | \ M \ | \ w \rangle$ will be defined as
\[
\langle v \ | \ M \ | \ w \rangle := \sum_{j,k=1}^d m_{jk}v_j w_k.
\]
In what follows, $\BDO (a;r)$ will denote the open ball in $\mathbb{R}^d$ centered at $a\in\mathbb{R}^d$
and having radius $r$, while $\BDC (a;r)$ will denote the closure of $\BDO (a;r)$. 
\medskip

Let $\rho$ be a defining function for $\bdy$. Recall that for $p\in\bdy$ and a vector $v \in T_p(\bdy)$, the 
second fundamental form for $\bdy$ at $p$ is the quadratic form
\[
T_p(\bdy) \ni v \mapsto \langle v \ | \ (\form\rho)(p) \ | \ v \rangle,
\]
where $\form\rho$ denotes the real Hessian of $\rho$.
We define $\Nu_p\subseteq T_p(\bdy)$ to be the null space of the second fundamental form at $p$, i.e. 
$\Nu_p = \{v\in T_p(\bdy) : \langle v \ | \ (\form\rho)(p) \ | \ v \rangle = 0\}$
\medskip

A final piece of notation : if $\phi$ is a $\mathcal{C}^1$ function defined in some open set in $\Cn$,
$\partial_k\phi$ and $\partial_{\bar k}\phi$ will denote
\[
\partial_k\phi = \frac{\partial\phi}{\partial z_k}, \qquad
		\partial_{\bar k}\phi = \frac{\partial\phi}{\partial\zbar_k}.
\]
\smallskip

A standard approach to proving that $\mn\subset\bdy$ is a peak-interpolation set --
$\mn$, $\bdy$ smooth and $\OM\Subset\Cn, \ n\geq 2$ -- which is encountered in the papers 
\cite{HT} and \cite{wR}, is to use Bishop's theorem \cite{eB}, which states :
\medskip

\begin{theorem}[Bishop]\label{T:Bishop} Let $\OM$ be a bounded domain in $\Cn$ and let $K\subset\bdy$ be a 
compact subset. If $K$ is a totally-null set -- i.e. if for every annihilating measure $\mu\perp A(\OM)$,
$|\mu|(K)=0$ -- then $K$ is a peak-interpolation set for $A(\OM)$.
\end{theorem}
\medskip

In the above theorem, an {\bf annihilating measure} refers to a regular, complex Borel measure on $\overline\OM$
which, viewed as a bounded linear functional on $\mathcal{C}(\overline\OM)$, annihilates $A(\OM)$.
\medskip

Bishop's theorem implies that it suffices to show that $\mn$ is a countable union of totally-null sets, which
is the approach taken in \cite{wR}. The essential difference between the proof of 
Theorem \ref{T:main} and the earlier results lies in the very particular manner in which we decompose 
$\mn\subset\bdy$, in the weakly convex setting, into countably many totally-null subsets. As we shall see,
the manner in which we decompose $\mn$ is necessitated by the fact that there may be submanifolds of $\mn$ 
along which the second fundamental form for $\bdy$ is {\em not} strictly positive -- a phenomenon that is absent
in the strictly convex setting.
\medskip

The proof of Theorem \ref{T:main} relies on four main ingredients. We need, for our proof, to show that :

\begin{enumerate}
\item[(1)] If $\OM\subset\RN, \ N\geq 2$, is a convex domain having a smooth boundary that contains no 
line segments, there cannot be a smooth curve $\sigma : I \to \bdy$ of class $\mathcal{C}^1$ with 
$\sigma^\prime(t)\in \Nu_{\sigma(t)}$ on an entire interval. Consequently -- as we will show in Section 3 -- if
$\mn\subset\bdy$ is a smooth submanifold, $\Nu_\zt\cap T_\zt(\mn)=\{0\}$ for each $\zt$ belonging to an 
open, dense subset of $\mn$.
\smallskip

\item[(2)] If $\OM\subset\RN, \ N\geq 2$, is a bounded convex domain with real-analytic boundary, 
$\mn\subset\bdy$ is a real-analytic submanifold and $p\in\mn$, there is a neighbourhood $V\ni p$ 
and a stratification of $\mn\cap V$ into finitely many real-analytic submanifolds (not necessarily closed)
of $\bdy\cap V$ such that if $\strat$ is a stratum of positive dimension, 
$T_\zt(\strat)\cap\Nu_\zt = \{0\} \ \forall\zt\in\strat$.
\smallskip

\item[(3)] For each stratum $\strat\subset\bdy$ of the aforementioned local stratification with 
$\Dim(\strat)\geq 1$, and for each $q\in\strat$, there is a small neighbourhood $U\ni q$ such that 
the compact $\strat\cap\overline{U}$ is a totally-null set. 
\end{enumerate}
\medskip

The central idea in \cite{wR} is to show that one can write $\mn = \cup_{j\in\mathbb{N}}K_j$, where
each $K_j$ is compact, in such a manner that each $K_j$ is totally-null. This relies on the 
ability to construct a family of functions $\{h_\delta\}_{\delta>0}\subset A(\OM)$ that is uniformly bounded on
$\overline\OM$, such that $h_\delta(z)\to 0$ as $\delta\to 0$, for each $z\in\OM$, and which, in the limit, has 
a specified behaviour on an $\mn$-open neighbourhood of $K_j$. The analogue of this construction, in our
context, is the following claim, which is valid in the more general setting of {\em smoothly} bounded, 
weakly convex domains. Item (3) above is a consequence of the following claim, which is the last key ingredient
in the proof of our main theorem.
\smallskip

\begin{enumerate}
\item[(4)] Let $\OM\subset\Cn, \ n\geq 2$, be a bounded, weakly convex domain having a smooth boundary 
that contains no line segments, and let $\gamma:\BDO(0;R)\to\bdy$ be a smooth imbedding 
whose image is complex-tangential. Also assume that $d\gamma(x)(\Rd)\cap\Nu_{\gamma(x)}=\{0\} \ \forall x$.
There exists a $\varrho>0$ such that if $f\in\mathcal{C}_{c}[\BDO(0;\varrho);\mathbb{C}]$, then defining
\[
h_\delta(z) = \int_{\BDO(0;\varrho)}\frac{\delta^d f(x)/G(x) \ dx}
				{\left\{\delta^2+\drv{\gamma(x)}{\gamma_j(x)-z_j}\right\}^d} \ , \ 
				z\in \overline\OM,
\]
where $G$ is defined as
\[
G(x) = \int_{\Rd}\left\{1+\frac{1}{2}\Q \right\}^{-d}dv \ ,
\]
we have : 
\begin{enumerate}
\item[$(i)$] $\{h_\delta\}_{\delta>0}\subset A(\OM)$ and is uniformly bounded on $\overline\OM$, 
\item[$(ii)$] $\lim_{\delta\to 0}h_\delta(z)=0$ if $z\in\overline\OM\setminus\gamma[\BDO(0;\varrho)]$,
\item[$(iii)$] $\lim_{\delta\to 0}h_\delta[\gamma(s)]=f(s) \ \forall s\in\BDO(0;\varrho)$.
\end{enumerate}
\end{enumerate}
\medskip

We remark that the object $\gamma[\BDC(0;\varrho)]$ -- $\gamma, \ \varrho$ as above -- is the prototye for 
compact sets of the sort described in item (3).  Furthermore, we observe that the family of integrals given 
above is the same as that appearing in \cite{wR}, although that paper is about a result similar to Theorem
\ref{T:main} but which applies only to {\em strictly} convex domains. Item (4) says that if for every 
$x\in\BDO(0;R)$ the second fundamental form for $\bdy$ is strictly positive on
$T_{\gamma(x)}\gamma[\BDO(0;R)]\subset T_{\gamma(x)}(\bdy)$, the aforementioned integrals have the estimates
similar to those in \cite{wR}. We present these estimates in Section 5, Theorem \ref{T:big}.
\medskip
  
In Section 3, we state and prove propositions relating to item (1) above. The geometric results from
Section 3 and real-analyticity are both crucial to the claim made in item (2) above. The existence of
a local stratification of $\mn$ having certain geometric properties is proved in Section 4 below. This 
local stratification is {\em essential} to our proof, and the need for it is what necessitates the hypothesis 
of real-analyticity in Theorem \ref{T:main}. We finally complete the proof of Theorem \ref{T:main} in Section 6.
\bigskip

\section{Results on Convex Domains}

\begin{lemma}\label{L:nonnull} Let $\OM$ be a convex domain in $\RN, \ N\geq 2$, having a
$\mathcal{C}^2$ boundary and 
containing no line segments in its boundary. Then, there is no smooth curve 
$\sigma : I \to \bdy$ of class $\mathcal{C}^1$ with $\sigma^\prime(t)\in \Nu_{\sigma(t)} \ \forall
t\in I$ (where $I$ is some interval of the real line).
\end{lemma}
\begin{proof}
Assume the result is false. Let $\sigma : I\to\bdy$ be a curve with 
$\sigma^\prime(t)\in \Nu_{\sigma(t)} \ \forall t\in I$ ($I$ is some interval).
Let $\rho$ be a defining function for $\bdy$ with $\|\nabla\rho\| = 1$.
For $t\in I$, set
\begin{align}
n(t) &= \nabla\rho(\sigma(t))\in \RN, \notag \\
H(t) &= (\form\rho)(\sigma(t))\in \mathbb{R}^{N\times N}. \notag
\end{align}
We compute that
\begin{equation}\label{E:chrule}
n^\prime(t) = H(t)\sigma^\prime(t).
\end{equation}
\smallskip

Now, notice that as $\sigma^\prime(t)\in \Nu_{\sigma(t)}$ by assumption, and as $H(t)$ is a symmetric matrix 
and is positive semi-definite on $T_{\sigma(t)}(\bdy) \ \forall t\in I$, we have
\[
\langle H(t)\sigma^\prime(t),v \rangle = \langle\sigma^\prime(t), \ H(t)v \rangle 
			= 0 \quad \forall v\in T_{\sigma(t)}(\bdy), \ \forall t\in I.
\]
The last equality  follows from the fact that, since $H(t)$ is positive semi-definite on 
$T_{\sigma(t)}(\bdy)$, for any  $v\in T_{\sigma(t)}(\bdy)$ we have 
\[
0\leq \langle\sigma^\prime(t)+\lambda v \ , \ H(t)(\sigma^\prime(t)+\lambda v) \rangle
	\ = \ \lambda^2\langle v, \ H(t)v \rangle + 2\lambda\langle\sigma^\prime(t), \ H(t)v \rangle 
	\quad \forall \lambda\in\mathbb{R},
\]
which forces $\langle\sigma^\prime(t), \ H(t)v \rangle$ to vanish. Thus, by \eqref{E:chrule}, $n^\prime(t)$ 
is orthogonal to $T_{\sigma(t)}(\bdy), \ \forall t\in I$.
\medskip

Next, observe that
\begin{align}
\langle n(t), n(t) \rangle &= 1 \notag \\
\Rightarrow \ 2\langle n^\prime(t), n(t) \rangle &= 0, 
			\quad\text{[by differentiating the above equation]} \notag
\end{align}
whence $n^\prime(t)$ is orthogonal to the outward unit normal at $\sigma(t)$ for
each $t\in I$. We infer, thus, that $n^\prime(t)=0 \ \forall t\in I$. Thus $n$ is 
constant on $I$.
\smallskip

Write $c = n(t)$, and define a function
\[
G(s,t) = \langle \sigma(s)-\sigma(t), c \rangle, \ s,t\in I.
\]
Clearly
\[
\frac{\partial G}{\partial s} = \frac{\partial G}{\partial t} = 0, 
				\quad\text{[since $\sigma^\prime(\centerdot)\perp n(\centerdot)$]}
\]
whence $G\equiv \text{const}$. Since $G(s,s)=0$, $G\equiv 0$. Thus,
\begin{equation}\label{E:line}
\sigma(s)-\sigma(t)\in T_{\sigma(t)}(\bdy).
\end{equation}
By the convexity of $\OM$, the line segment joining $\sigma(s)$ and $\sigma(t)$ must
lie in $\overline{\OM}$. In view of \eqref{E:line}, this means that the line
segment joining $\sigma(s)$ to $\sigma(t)$ lies in $\bdy$. This is a contradiction,
whence the initial assumption if false.
\end{proof}

\begin{lemma}\label{L:density} Let $\OM$ be a convex domain in $\RN, \ N\geq 2$, having a $\mathcal{C}^2$
boundary and containing no line segments in its boundary, and let $\mn$ be a submanifold of $\bdy$
of class $\mathcal{C}^2$. Then, the set $\{p\in\mn \ | \ T_p(\mn)\cap \Nu_p = \{0\} \ \}$
is open and dense in $\mn$.
\end{lemma}
\begin{proof}
Let $\text{dim}_{\mathbb{R}}(\mn)=d>0$. Define
\[
\mathfrak{S} = \{p\in\mn \ | \ T_p(\mn)\cap \Nu_p \varsupsetneq \{0\} \ \}.
\]
Let $\gamma : (\BDO(0;\eps),0)\to (\mn,p)$ be a non-singular parametrization of $\mn$ near $p\in\mn$ 
of class $\mathcal{C}^2$. We will show that $\mathfrak{S}\cap\gamma[\BDO(0;\eps)]$ cannot contain an 
open subset of $\mn$. 
Define
\begin{align}
H(s) &= d\gamma(s)^{T}(\form\rho)(\gamma(s)) \ d\gamma(s), \notag \\
\mathcal{N}_s &= \{v\in\mathbb{R}^d : \langle v \ | \ H(s) \ | \ v \rangle = 0\}, \notag
\end{align}
Consider the function $G : \BDO(0;\eps)\times \Rd\setminus\{0\}\to \mathbb{R}$ defined by
\[
G : (s \ ; \ v_1,...,v_d) \mapsto 
	\langle d\gamma(s)v \ | \ (\form\rho)(\gamma(s)) \ | \ d\gamma(s)v \rangle .
\]
Observe that the matrix of derivatives with respect to $v$ 
\[
d_vG(s,v) = 2v^{T}H(s)\in \mathbb{R}^{1\times d}.
\]
\smallskip

Assume that $\text{int}[\mathfrak{S}\cap\gamma[\BDO(0;\eps)]]\neq\emptyset$. Without loss of generality,
we may assume that there exists an $\eps_*\in (0,\eps]$ such that 
$\gamma[\BDO(0;\eps_*)]\subset\mathfrak{S}\cap\gamma[\BDO(0;\eps)]$. We then have the following
situations :
\medskip

\noindent{{\em Case (i).} There is a $v^0\in S^{d-1}$ and $s_0\in\BDO(0;\eps_*)$ such that
$d_vG(s_0,v^0)\neq 0$. This means that $\Rnk G(s_0,v^0)$ is maximal. Therefore, by the implicit function
theorem, there is a $\delta\in(0,\eps_*]$ and a non-vanishing vector-field 
$F=(F_1,...,F_d) : \BDO(s_0;\delta) \to \Rd\setminus\{0\}$ of class $\mathcal{C}^2$ such that
\begin{align}
G(s,F(s)) &= 0 \quad \forall s\in\BDO(s_0;\delta), \notag \\
F(s_0) &= v^0. \notag
\end{align}}
\smallskip

\noindent{{\em Case (ii).} For each $s\in\BDO(0;\eps_*)$ and $v\in\Rd\setminus\{0\}$, $d_vG(s,v)=0$.
In this situation, $\mathcal{N}_s=\Rd \ \forall s\in\BDO(0;\eps_*)$. In this case let 
$F=(F_1,...,F_d) : \BDO(0;\eps_*) \to \Rd\setminus\{0\}$ be any non-vanishing vector field of class
$\mathcal{C}^2$.}
\medskip

In either of the above cases, we have
\begin{equation}\label{E:field}
F_1(s)\frac{\partial\gamma}{\partial s_1}(s)+\dots +
	F_d(s)\frac{\partial\gamma}{\partial s_d}(s) \in \Nu_{\gamma(s)}
\quad \forall s\in B(s_0;\delta).
\end{equation}
Let $\sigma : (-a,a)\to \BDO(s_0;\delta)$ be the integral curve to $F$ through $s_0$, i.e.
\begin{align}
\sigma(0) &= \ s_0, \notag\\
\sigma^\prime(t) &= \ F(\sigma(t)) \quad \forall t\in (-a,a). \notag
\end{align}
Then
\begin{equation}\label{E:key}
(\gamma\circ\sigma)^\prime (t)
	= F_1(\sigma(t))\frac{\partial\gamma}{\partial s_1}(\sigma(t))+\dots +
F_d(\sigma(t))\frac{\partial\gamma}{\partial s_d}(\sigma(t)) \quad \forall t\in (-a,a).
\end{equation}
From \eqref{E:field} and \eqref{E:key},
\[
(\gamma\circ\sigma)^\prime (t) \in \Nu_{\gamma\circ\sigma(t)} \quad \forall t\in (-a,a)
\]
which is impossible, by Lemma \ref{L:nonnull}. Thus, $\mathfrak{S}$ does not contain any open subsets 
of $\mn$.
\medskip

In particular $\mn\setminus\mathfrak{S}\neq\emptyset$. Consider a point $p\in\mn\setminus\mathfrak{S}$,
and let $\gamma : (\BDO(0;\eps),0)\to(\mn,p)$ be as before.
Consider $G : \BDO(0;\eps)\times S^{d-1}\to \mathbb{R}$ defined exactly as above.
$G^{-1}[\mathbb{R}\setminus\{0\}]$ is an open set and 
$G^{-1}[\mathbb{R}\setminus\{0\}]\supset \{0\}\times S^{d-1}$, since $p\in\mn\setminus\mathfrak{S}$. 
From this, we infer that there is an $\mn$-open neighbourhood of $p$ contained in $\mn\setminus\mathfrak{S}$.
This last fact completes the proof.
\end{proof}
\smallskip

\section{A Stratification Theorem}

In this section, we shall state precisely, and prove, the informally stated fact in item (2) in Section 2. A
key fact that we will use to do so is the structure theorem for real-analytic subvarieties of $\RN, \ N\geq 2$.
This theorem is due to {\L }ojasiewicz \cite{sL}, which we now state.

\begin{theorem}[{\L }ojasiewicz]\label{T:Loja} Let $F$ be a non-constant real-analytic function defined in a
neighbourhood of $0\in\RN$, and assume that $V(F)=F^{-1}\{0\}\ni 0$. Then, there is a 
small neighbourhood $U\ni 0$ such that $V(F)\cap U$ has the decomposition
\[ 
V(F)\cap U = \cup_{j=0}^{N-1} S_j,
\] 
where each $S_j$ is a finite, disjoint union of (not necessarily closed) $j$-dimensional real-analytic 
submanifolds contained in $U$, such that each connected component of $S_j$ is a closed real-analytic submanifold 
of $U\setminus \left(\cup_{k=0}^{j-1}S_k\right), \ j=1,...,(N-1)$.
\end{theorem}
\smallskip

We remark that although the above theorem describes the local structure of the zero-set of a single 
real-analytic function, it, in fact, describes the local structure of a variety near the origin. This is
because, given finitely many real-analytic functions $f_1,...,f_M$ that vanish at the origin, their set
of common zeros is exactly the zero-set of the real-analytic function $F:= |f_1|^2+\dots +|f_M|^2$.
\medskip

The following theorem is a precise statement of item (2) in Section 2.
\medskip

\begin{theorem}\label{T:key} Let $\OM$ be a bounded convex domain in $\RN, \ N\geq 2$, having real-analytic
boundary, and let $\mn\subset\bdy$ be a $d$-dimensional real-analytic submanifold. Let $p\in\mn$. There is
a neighbourhood $V\ni p$ such that
\begin{equation}\label{E:strata}
\mn\cap V = \cup_{j=0}^d M_j,
\end{equation}
where
\begin{enumerate}
\item[(i)] Each $M_j$ is a disjoint union of finitely many (not necesarily closed) 
$j$-dimensional real-analytic submanifolds contained in $V$.
\item[(ii)] Each connected component of $M_j$ is a closed, real-analytic submanifold of 
$V\setminus\left(\cup_{k=0}^{j-1}M_k\right), \ j=1,...,d$.
\item[(iii)] For any $j\neq 0$, if $M_{j,\alpha}$ is a connected component of $M_j$, 
$\Nu_\zt\cap T_\zt(M_{j,\alpha}) = \{0\} \ \forall \zt\in M_{j,\alpha}$.
\end{enumerate}
\end{theorem}
\begin{proof}
Fix $p\in\mn$. Let $\gamma : (\BDO(0;\eps),0)\to (\mn, p)$ be a real-analytic parametrization 
of $\mn$ near $p$ such that $\Rnk[d\gamma(x)]$ is maximal $\forall x$. Consider the real-analytic 
function $\mathfrak{F} : \BDO(0;\eps)\to\mathbb{R}$ defined by
\[
\mathfrak{F}(x)
= \text{det}\left[d\gamma(x)^{T}(\form\rho)(\gamma(x)) \ d\gamma(x)\right].
\]
The matrix in the above expression is simply the pull-back of the Hessian $\form\rho$ by $\gamma$.
By Lemma \ref{L:density}, $\mathfrak{F}\not\equiv 0$. Without loss of generality, we may assume that
if $\mathfrak{F}^{-1}\{0\}\neq\emptyset$, then $\mathfrak{F}^{-1}\{0\}\ni 0$. By {\L }ojasiewicz's theorem 
\cite{sL}, there is a neighbourhood $U\ni 0$, $U\subseteq\BDO(0;\eps)$, such that
\begin{equation}\label{E:strat1}
\mathfrak{F}^{-1}\{0\}\cap U = \cup_{j=0}^{d-1}\STrat_j,
\end{equation}
where each $\STrat_j$ is a disjoint union of finitely many $j$-dimensional real-analytic submanifolds,
and each connected component of $\STrat_j$ is a closed submanifold of  
$U\setminus\left(\cup_{k=0}^{j-1}\STrat_k\right), \ j=1,...,d-1$. Write 
$\STrat_d = U\setminus\left(\cup_{j=0}^{d-1}\STrat_j\right)$.
\medskip

We plan to demonstrate the present result by induction. We make the following inductive hypothesis :
\smallskip

For $m<d-1$ we have, shrinking $U$ if necessary, a stratification of $U$
\begin{equation}\label{E:strat2}
U = \cup_{j=0}^d S_j,
\end{equation} 
where,
\begin{enumerate}
\item[$(a) \ \text{ }$] Each $S_j$ is a disjoint union of finitely many (not necesarily closed) 
$j$-dimensional real-analytic submanifolds contained in $U$.
\item[$(b) \ \text{ }$] Each connected component of $S_j$ is a closed, real-analytic submanifold of 
$U\setminus \left(\cup_{k=0}^{j-1}S_k\right), \ j=1,...,d$.
\item[$(c)_m$] For each $k=0,...,m$ and each connected component $S_{d-k,\alpha}$, of $S_{d-k} \ $,
$\tansp{d-k}{\gamma(x)}\cap\Nu_{\gamma(x)} = \{0\} \ \forall x\in S_{d-k,\alpha}$.
\item[$(d)_m$] $\cup_{j=0}^{(d-m-1)}S_j$ is a real-analytic subvariety of $U$.
\end{enumerate}
\smallskip

Consider the real-analytic subvariety $\widetilde{V}$ of $U$ given by
\[
\widetilde{V} = \left(\cup_{j=0}^{(d-m-1)}S_j\right)\cap
		\{x\in U : \Rnk\left[d\gamma(x)^{T}(\form\rho)(\gamma(x)) \ d\gamma(x)\right]
		\leq(d-m-2)\},
\]
where the $S_j$'s come from the stratification in \eqref{E:strat2}.
We consider $S_{d-m-1,\alpha}$ : a connected component of $S_{d-m-1}$. 
$\umm = \gamma(S_{d-m-1,\alpha})$ is a real-analytic submanifold contained in $\bdy$.
By Lemma \ref{L:density}, there is an $\umm$-open set $\mathcal{V}$ such that 
$T_\zt(\umm)\cap\Nu_\zt = \{0\}, \ \forall \zt\in\mathcal{V}$. 
Write $\mathcal{U}= \left(\gamma |_{S_{d-m-1,\alpha}}\right)^{-1}(\mathcal{V})$. $\mathcal{U}$ is
open in $S_{d-m-1,\alpha}$, and for any $x\in \mathcal{U}$ and any 
$v\in [T_x(S_{d-m-1,\alpha})\setminus\{0\}] \ $,
$d\gamma(x)v \notin \Nu_{\gamma(x)}$. In other words, for each $x\in\mathcal{U}$,
\[
\text{ker}\left\{d\gamma(x)^{T}(\form\rho)(\gamma(x)) \ d\gamma(x)|_{T_x(S_{d-m-1,\alpha})}\right\}=\{0\},
\]
where we identify the matrices $d\gamma(x)^{T}(\form\rho)(\gamma(x)) \ d\gamma(x)$ with linear
transformations. So, for each $S_{d-m-1,\alpha}$, the real-analytic subvariety 
$(\widetilde{V}\cap S_{d-m-1,\alpha})\subsetneq S_{d-m-1,\alpha}$. From this, we infer that 
$\Dim(\widetilde{V})<(d-m-1)$. By {\L }ojasiewicz's theorem, shrinking $U$ if necessary, we have
\[
\widetilde{V}\cap U = \cup_{j=0}^{(d-m-2)}\varST_j
\]
where each connected component of $\varST_j$ is a closed submanifold of 
$U\setminus\left(\cup_{k=0}^{j-1}\varST_k\right), \ j=1,...,d-m-2$. Now write
\begin{equation*}
\st_j = \begin{cases} S_j\cap U, & \text{if $j\geq(d-m)$}\\
			(S_{d-m-1,\alpha}\setminus\widetilde{V})\cap U, & \text{if $j=(d-m-1)$}\\
			(S_j\cup\varST_j)\cap U, & \text{if $j\leq(d-m-2)$},
	\end{cases}
\end{equation*}
shrinking $U$ further if necessary so that 
\[
U = \cup_{j=0}^d \st_j
\]
is a stratification of $U$ that satisfies $(a)$ and $(b)$ above with $\st_j$ replacing $S_j$.
By construction, each connected component $\st_{d-m-1,\alpha}$, of $\st_{d-m-1}$, satisfies
$(c)_{m+1}$ and $\left(\cup_{j=0}^{(d-m-2)}\st_j\right)= \widetilde{V}\cap U$ satisfies $(d)_{m+1}$.
\smallskip

Notice that the stratification in \eqref{E:strat1} establishes the case $m=0$ for the inductive hypothesis
above. By induction, therefore, we can find, shrinking $U$ if necessary, a stratification
\begin{equation}\label{E:strat3}
U=\cup_{j=0}^d S_j
\end{equation}
where each connected component $S_{j,\alpha}$, of $S_j, \ j=1,...,d$, is a closed, real-analytic
submanifold of $U\setminus \left(\cup_{k=0}^{j-1}S_k\right)$, and for each $j\geq 1$ and each $\alpha$,
$\tansp{j}{\zt}\cap\Nu_\zt=\{0\}, \ \forall \zt\in\gamma(S_{j,\alpha})$. We now find a suitably small
neighbourhood, say $V$, of $p$ so that writing
\[
M_j = \gamma(S_j)\cap V, \qquad \mn\cap V = \cup_{j=0}^d M_j
\]
[where the $S_j$'s come from \eqref{E:strat3}] gives us the result.
\end{proof}
\medskip

\section{Quantitative Results}

In this section, we work with bounded convex domains $\OM\subset\Cn, \ n\geq 2$, having smooth 
boundaries containing no line segments. Let $\gamma:(\BDO(0;R),0)\to(\bdy,q)$ be a smooth imbedding 
whose image is complex-tangential, and for which $d\gamma(x)(\Rd)\cap\Nu_{\gamma(x)}=\{0\} \ \forall x$.
For the remainder of this section, $\gamma$ and $R>0$ will have the specific meaning just introduced.
In the context of Theorem \ref{T:main}, given a point $p\in\mn$, $\gamma[\BDO(0;R)]$ serves as the
prototype for an open subset of a stratum of positive dimension in the local stratification 
\eqref{E:strata} of $\mn$ near $p$.
\medskip

For $\OM$ as above, $\rho$ a smooth defining function for $\bdy$, $\zt\in\bdy$ and $z\in\Cn$, we
write $G(\zt,z) = \vardrv{\zt}{\zt_j-z_j}$. For a fixed $\zt\in\bdy$, the equation $G(\zt,z)=0$
defines $H_\zt(\bdy)$, and the real part of $G(z,\zt)$ is the perpendicular distance of $z$ from 
$T_\zt(\bdy)$. Thus, by the convexity of $\OM$, if $z\in\overline\OM$, $\er[G(\zt,z)]\geq0$ and 
$G(\zt,z)=0 \ \Leftrightarrow \ z=\zt$. In other words, $\{G(\zt,\centerdot)\}_{\zt\in\bdy}$ is a 
smoothly varying family of peak functions for $A(\OM)$.
\medskip

We now prove a technical lemma, which we will need later in this section.
\medskip

\begin{lemma}\label{L:localmin} Let $\OM$, $\gamma$ and $R$ be as described above. For each 
$r\in(0,R/2)$, there exists an open set $\mathcal{U}(r)\supset\gamma[\BDO(0;2r)]$ such that for each
$z\in\overline\OM\cap\mathcal{U}(r)$, there exists a $\yz\in\BDO(0;2r)$ satisfying
\begin{multline}
\er\left\{\drv{\gamma(\yz)}{\gamma_j(\yz)-z_j}\right\}
	\leq \er\left\{\drv{\gamma(x)}{\gamma_j(x)-z_j}\right\}\quad  \text{for every $x$ belonging} \\
\text{to a small neighbourhood, $U_z\subseteq\BDO(0;R)$, of $\yz$.}
\end{multline}
\end{lemma}
\begin{proof}
In what follows, we will write $\er\left\{\vardrv{\zt}{\zt_j-z_j}\right\}=F(\zt,z)$. For $\zt\in
\gamma[\BDO(0;R)]$, we define
\[
N_\zt(\gamma;\eps) := \{z\in N_\zt(\gamma[\BDO(0;R)]) : |z-\zt|<\eps\},
\]
where $N_\zt(\gamma[\BDO(0;R)])$ denotes the normal space of $\gamma[\BDO(0;R)]$ in $\Cn$ at $\zt$. 
Let $\sigma>0$ be so small that if $z$ lies in a tube around $\gamma[\BDO(0;R)]$, 
\[
\text{dist}[z,\gamma[\BDO(0;R)]]\leq\sigma \ \Rightarrow \ \text{there is a unique $x\in\BDO(0;R)$ 
such that $z\in\overline{N_{\gamma(x)}(\gamma;\sigma)}$.}
\]
Also we will assume (shrinking $R>0$ if necessary) that for each $r\in(0,R/2)$ and each $x\in\BDO(0;2r)$,
$N_{\gamma(x)}(\gamma;2\sigma)\cap\gamma[\partial\BDC(0;2r)]=\emptyset$.
We now fix $r\in(0,R/2)$ for the remainder of this proof. For each $t\in (0,2r)$, define the 
function 
$\mathfrak{F}_t : \overline\OM\cap\left(\cup_{|x|\leq t}\overline{N_{\gamma(x)}(\gamma;\sigma)}\right)
\times\gamma[\partial\BDC(0;2r)]
\to [0,\infty)$ by
\[
\mathfrak{F}_t : (z,\xi)\mapsto F(\xi,z).
\]
For a fixed $t\in (0,2r)$, $\mathfrak{F}_t(z,\xi)>0$, by convexity and by the foregoing choices 
for $R$ and $\sigma$. Thus, there
exists a $m_t>0$ such that $\mathfrak{F}_t(z,\xi)\geq m_t \ \forall (z,\xi)\in
\overline\OM\cap\left(\cup_{|x|\leq t}\overline{N_{\gamma(x)}(\gamma;\sigma)}\right)
\times\gamma[\partial\BDC(0;2r)]$. 
Write $s_t=\min\{m_t/2,\sigma\}$. Then
\begin{multline}\label{E:estimate}
z\in \overline\OM\cap\left(\cup_{|x|\leq t}\overline{N_{\gamma(x)}(\gamma;s_t)}\right) \\
\Rightarrow \ F(\zt^z,z)=\text{dist}[z,T_{\zt^z}(\bdy)]\leq s_t < F(\xi,z), \quad
\forall \xi\in\gamma[\partial\BDC(0;2r)],
\end{multline}
where $\zt^z\in\gamma[\BDC(0;t)]$ such that $z\in N_{\zt^z}(\gamma[\BDO(0;R)])$.
We define
\[
\mathcal{U}(r) = \text{int}\left[\cup_{t\in(0,2r)}\left\{\cup_{|x|\leq t}N_{\gamma(x)}(\gamma;s_t)\right\}\right],
\]
and for each $z\in\overline\OM\cap\mathcal{U}(r)$, we define $\yz$ by
\[
F(\gamma(\yz),z) = \inf_{x\in\BDC(0;2r)}F(\gamma(x),z).
\]
If we could show that $\yz\notin\partial\BDC(0;2r)$, then we would be done. For each 
$z\in\mathcal{U}(r)$, let $x^z\in\BDO(0;2r)$ be such that $z\in N_{\gamma(x^z)}(\gamma[\BDO(0;R)])$. If 
$z\in\cup_{|x|\leq t}N_{\gamma(x)}(\gamma;s_t)$ for some $t\in (0,2r)$, then $|x^z|\leq t$. In view
of \eqref{E:estimate}
\[
F(\gamma(\yz),z) \ \leq \ F(\gamma(x^z),z) \ < \ F(\gamma(s),z) \quad 
						\forall s\in\partial\BDC(0;2r).
\]
Hence, $\yz\notin\partial\BDC(0;2r)$, and we have our result.
\end{proof}
\medskip

So far, we have not made use of the fact that $\gamma[\BDO(0;R)]$ is complex-tangential. We shall do so in 
the next three lemmas.
\medskip

\begin{lemma}\label{L:critical} Let $\OM$, $\gamma$, $\yz\in\BDO(0;2r)$ and $\mathcal{U}(r)$ be as 
in Lemma \ref{L:localmin}. Then, for $z\in\overline\OM\cap\mathcal{U}(r)$
\begin{equation}\label{E:vanishing}
\er\left\{\adda\dbdrv[\gamma(\yz)]\gamder{k}{\mu}(\yz)[\gamma_j(\yz)-z_j] + 
		\dbdrvc[\gamma(\yz)]\overline{\gamder{k}{\mu}(\yz)}[\gamma_j(\yz)-z_j]\right\}
		= 0, \ \mu=1,...,d.
\end{equation}
\end{lemma}
\begin{proof}
If $z\in\overline\OM\cap\mathcal{U}(r)$, then $\yz$ is a local minimum of the function
\[
\BDO(0;2r)\ni x\mapsto F(\gamma(x),z).
\]
Therefore, taking the partial derivative of the above with respect to $x_\mu$ and evaluating at
$x=\yz$, we get
\begin{multline}\label{E:partial}
\er\left\{\adda\dbdrv[\gamma(\yz)]\gamder{k}{\mu}(\yz)[\gamma_j(\yz)-z_j] + 
		\adda\dbdrvc[\gamma(\yz)]\overline{\gamder{k}{\mu}(\yz)}[\gamma_j(\yz)-z_j]\right. \\
		\left. + \sum_{j=1}^n\partial_j\rho(\gamma(\yz))\gamder{j}{\mu}(\yz)\right\}
				= 0, \ \mu=1,...,d.
\end{multline}
Since $\gamma[\BDO(0;R)]$ is complex-tangential, we have
\begin{equation}\label{E:comptang}
\sum_{j=1}^n\partial_j\rho(\gamma(\yz)) \ \gamder{j}{\mu}(\yz) = 0, \ \mu=1,...,d.
\end{equation}
The result follows from \eqref{E:partial} and \eqref{E:comptang}.
\end{proof}
\smallskip

In the next lemma, we exploit the fact that $d\gamma(x)(\Rd)\cap\Nu_{\gamma(x)}=\{0\} \ \forall x\in\BDO(0;R)$.
\medskip

\begin{lemma}\label{L:pos} Let $\OM$ and $\gamma$ be as described above. There 
exist uniform constants $\delta\equiv\delta(\gamma)>0$ and $C\equiv C(\gamma)>0$ such that
\begin{equation}\label{E:positivity}
\er\left\{\drv{\gamma(x)}{\gamma_j(x)-\gamma_j(y)}\right\} \geq C \ |x-y|^2 \ \forall
									x,y\in\BDO(0,\delta).
\end{equation}
\end{lemma}
\begin{proof}
Let $\eta\in\bdy$ and $z\in\overline\OM$. We Taylor expand the function 
\[
\eta \mapsto \vardrv{\eta}{\eta_j-z_j}
\]
about $\eta=\xi$ to get
\begin{align}
\vardrv{\eta}{\eta_j-z_j} &= \vardrv{\xi}{\xi_j-z_j} + \vardrv{\xi}{\eta_j-\xi_j} + 
				\adda\dbdrv(\xi)(\xi_j-z_j)(\eta_k-\xi_k) \notag \\
	&\quad + \adda\dbdrvc(\xi)(\xi_j-z_j)(\overline{\eta_k}-\overline{\xi_k}) 
			+ \adda\dbdrv(\xi)(\eta_j-\xi_j)(\eta_k-\xi_k) \notag \\
	&\quad + \adda\dbdrvc(\xi)(\eta_j-\xi_j)(\overline{\eta_k}-\overline{\xi_k}) +
					O(|\xi-z||\eta-\xi|^2, \ |\eta-\xi|^3). \notag
\end{align}
Substituting $z=\gamma(y), \ \eta=\gamma(x) \ \text{and} \ \xi=\gamma(y)$ in the above expression,
we get
\begin{multline}\label{E:expand1}
\drv{\gamma(x)}{\gamma_j(x)-\gamma_j(y)} = \drv{\gamma(y)}{\gamma_j(x)-\gamma_j(y)} +
		\adda\dbdrv(\gamma(y))[\gamma_j(x)-\gamma_j(y)][\gamma_k(x)-\gamma_k(y)]  \\
+ \adda\dbdrvc(\gamma(y))[\gamma_j(x)-\gamma_j(y)][\overline{\gamma_k(x)}-\overline{\gamma_k(y)}] +
		O(|\gamma(x)-\gamma(y)|^3),
\end{multline}
for $x,y\in\BDO(0;R)$. Taylor expanding $\gamma$ around $y\in\BDO(0;R)$, and substituting in
\eqref{E:expand1}, we have
\begin{align}\label{E:expand2}
\drv{\gamma(x)}{\gamma_j(x)-\gamma_j(y)} &= 
		\add\partial_j\rho(\gamma(y))\left\{\zum\gamder{j}{\mu}(y)(x_\mu-y_\mu)+
		\dfrac{1}{2}\zuma\fuder(x_\mu-y_\mu)(x_\nu-y_\nu)\right\} \\
	&\quad + \adda\dbdrv(\gamma(y))\left\{\zuma\gamder{j}{\mu}(y) \ \gamder{k}{\nu}(y)
							(x_\mu-y_\mu)(x_\nu-y_\nu)\right\} \notag \\
	&\quad + \adda\dbdrvc(\gamma(y))\left\{\zuma\gamder{j}{\mu}(y) \ \overline{\gamder{k}{\nu}(y)}
							(x_\mu-y_\mu)(x_\nu-y_\nu)\right\} \notag \\
	&\quad + O(|x-y|^3) \quad \forall x,y\in\BDO(0;R). \notag
\end{align}
\smallskip

We now use the fact that $\gamma[\BDO(0;R)]$ is complex-tangential. For $y\in\BDO(0;R)$, we have
\begin{equation}\label{E:comptang*}
\add\partial_j\rho(\gamma(y)) \ \gamder{j}{\mu}(y)=0, \quad \mu=1,...,d.
\end{equation}
Differentiating the above expression with respect to $x_\nu$ and evaluating at $x=y$ yields
\begin{multline}\label{E:comptang**}
\adda\dbdrv(\gamma(y))\gamder{j}{\mu}(y) \ \gamder{k}{\nu}(y)+\adda\dbdrvc(\gamma(y))
					\gamder{j}{\mu}(y) \ \overline{\gamder{k}{\nu}(y)} \\
+ \add\partial_j\rho(\gamma(y))\fuder = 0, \quad \mu,\nu = 1,...,d. 
\end{multline}
From \eqref{E:expand2}, \eqref{E:comptang*} and \eqref{E:comptang**}, we get
\[
\er\left\{\drv{\gamma(x)}{\gamma_j(x)-\gamma_j(y)}\right\} = 
	\frac{1}{2}\langle \ d\gamma(y)(x-y) \ | \ (\form\rho)(\gamma(y)) \ | \ d\gamma(y)(x-y) \ \rangle
			+ O(|x-y|^3).
\]
This statement, in conjunction with the strict positivity of $(\form\rho)(\zt)$ on 
$T_\zt(\gamma[\BDO(0;R)])\subset T_\zt(\bdy) \ \forall \zt\in\gamma[\BDO(0;R)]$, allows us to infer that there 
are uniform constants
$\delta\equiv\delta(\gamma)>0$ and $C\equiv C(\gamma)>0$ such that
\[
\er\left\{\drv{\gamma(x)}{\gamma_j(x)-\gamma_j(y)}\right\} \geq C \ |x-y|^2
						\quad \forall x,y\in\BDO(0;\delta).
\]
\end{proof}
\smallskip

\begin{lemma} Let $\OM$, $\gamma$ and $R$ be as described above. 
Then
\begin{multline}\label{E:keylimit}
\lim_{\delta\to 0} \ \er\left\{\dfrac{1}{\delta^2}\drv{\gamma(x+\delta v)}
		{\gamma_j(x+\delta v)-\gamma_j(x)}
		\right\} \\
= \dfrac{1}{2}\Q
\end{multline}
for any $x\in\BDO(0;R)$ and $v\in\mathbb{R}^d$.
\end{lemma}
\begin{proof}
Follows from \eqref{E:expand2}, \eqref{E:comptang*} and \eqref{E:comptang**} in the proof of Lemma 
\ref{L:pos}
\end{proof}
\medskip

\begin{theorem}\label{T:big} Let $\OM\subset\Cn, \ n\geq 2$, be a bounded, weakly convex domain having a 
smooth boundary that contains no line segments, and let $\gamma:\BDO(0;R)\to\bdy$ be a smooth imbedding 
whose image is complex-tangential. Also assume that $d\gamma(x)(\Rd)\cap\Nu_{\gamma(x)}=\{0\} \ \forall x$.
There exists a $\varrho\equiv\varrho(\gamma)>0$ such that if $f\in\mathcal{C}_{c}[\BDO(0;\varrho);\mathbb{C}]$, 
then defining
\begin{equation}\label{E:integ1}
h_\delta(z) = \int_{\BDO(0;\varrho)}\frac{\delta^d f(x)/G(x) \ dx}
				{\left\{\delta^2+\drv{\gamma(x)}{\gamma_j(x)-z_j}\right\}^d} \ , \ 
				z\in \overline\OM,
\end{equation}
where $G$ is defined as
\begin{equation}\label{E:integ2}
G(x) = \int_{\Rd}\left\{1+\frac{1}{2}\Q \right\}^{-d}dv \ ,
\end{equation}
we have
\begin{enumerate}
\item[$(i)$] $\{h_\delta\}_{\delta>0}\subset A(\OM)$ and is uniformly bounded on $\overline\OM$,
\item[$(ii)$] $\lim_{\delta\to 0}h_\delta(z)=0$ if $z\in\overline\OM\setminus\gamma[\BDO(0;\varrho)]$,
\item[$(iii)$] $\lim_{\delta\to 0}h_\delta[\gamma(s)]=f(s) \ \forall s\in\BDO(0;\varrho)$.
\end{enumerate}
\end{theorem}
\begin{proof} Consider any $r\in(0,R/3)$ and any $z\in\overline\OM\cap\mathcal{U}(r)$, 
where $\mathcal{U}(r)$ is as described in Lemma \ref{L:localmin}. We first estimate the quantity 
\[
\er\left\{\add\left[\derrho{\gamma(\yz+x)}-\derrho{\gamma(\yz)}\right][\gamma_j(\yz)-z_j]\right\},
\]
(where $\yz$ is as introduced in Lemma \ref{L:localmin}) given that $(\yz+x)\in\BDO(0;R)$.
Taylor expanding about $\yz$ and using the complex-tangency hypothesis for $\gamma[\BDO(0;R)]$, we have
\begin{multline}
\er\left\{\add\left[\derrho{\gamma(\yz+x)}-\derrho{\gamma(\yz)}\right][\gamma_j(\yz)-z_j]\right\} \\
	= \er\left\{\zum\left[\adda\dbdrv[\gamma(\yz)]\gamder{k}{\mu}(\yz)[\gamma_j(\yz)-z_j] + 
		\dbdrvc[\gamma(\yz)]\overline{\gamder{k}{\mu}(\yz)}[\gamma_j(\yz)-z_j]\right]
		x_\mu\right\} + O(|x|^2)|\gamma(\yz)-z|. \notag 
\end{multline}
In view of Lemma \ref{L:critical}, we can find uniform constants $c>0$ and $\eps_*>0$ such that
\begin{align}\label{E:secOrdEst*}
\left|\er\left\{\add\left[\derrho{\gamma(\yz+x)}-\derrho{\gamma(\yz)}\right][\gamma_j(\yz)-z_j]\right\}\right|
	\ \leq \ c|x|^2 \ |\gamma(\yz)-z|	& \quad \forall r\in(0,R/3), \\
					 { }	& \quad \forall |x|\leq\eps_*, \notag \\
					 { }	& \quad \forall z\in\overline\OM\cap\mathcal{U}(r). \notag
\end{align}
Let $\rod\in(0,R/3)$ be so small that for every $r\in(0,\rod]$,
\begin{align}
|\gamma(x)-z|\leq\frac{C(\gamma)}{2c}	& \quad \forall x\in\BDO(0;r), \notag \\
				 { }	& \quad \forall z\in\overline\OM\cap\mathcal{U}(r), \notag
\end{align}
where $c$ is as in \eqref{E:secOrdEst*} and $C(\gamma)$ is the constant appearing in Lemma \ref{L:pos}.
We now define a constant
\[
\varrho\equiv\varrho(\gamma):=\min\{1/2,\rod,\delta(\gamma)/2,\eps_*/3\},
\] 
where $\delta(\gamma)$ is the constant appearing in Lemma \ref{L:pos}. In what follows, we will use
the notation $y^z\in\BDO(0;2\varrho)$ to mean $y^z=y^z_{\varrho}$. From the preceding estimate and 
\eqref{E:secOrdEst*}, we have the estimate
\begin{align}\label{E:secOrdEst}
\left|\er\left\{\add\left[\derrho{\gamma(y^z+x)}-\derrho{\gamma(y^z)}\right][\gamma_j(y^z)-z_j]\right\}\right|
	\ \leq \ \frac{C(\gamma)}{2}\ |x|^2	& \quad \forall |x|\leq\eps_*, \\
					 { }	& \quad \forall z\in\overline\OM\cap\mathcal{U}(\varrho).\notag
\end{align}
\smallskip

Now, consider any $f\in\mathcal{C}_c[\BDO(0;\varrho);\mathbb{C}]$. For each $\delta>0$, define
$h_\delta$ according to \eqref{E:integ1} above. We remark that by our assumption on $\gamma$, the form
\[
\Rd\ni v\mapsto \er \Q
\]
is strictly positive definite for each $x\in\BDO(0;\varrho)$. Using this fact, it can be shown --
see \cite[Lemma 2.4]{wR} --that the integrals \eqref{E:integ2} converge, and that $G(x)\neq 0$.
From the discussion at the beginning of this section, we conclude that
the real part of $ \ \drv{\gamma(x)}{\gamma_j(x)-z_j}$, which occurs in the denominator
of the integral in \eqref{E:integ1} is non-negative when $z\in\overline\OM$. Thus, 
$h_\delta\in A(\OM)$ for each $\delta>0$.
\medskip

\noindent{{\bf Claim (i)} $\{h_\delta\}_{\delta>0}$ is uniformly bounded on $\overline\OM$.}
\smallskip

\noindent{We first consider the case when $z\in\overline\OM\cap\mathcal{U}(\varrho)$. We indulge in a 
slight abuse of notation : we will define the integrand in \eqref{E:integ1} to be $0$ when 
$x\notin\BDO(0;\varrho)$, whence $\BDO(0;\varrho)$ may be replaced by $\Rd$ in \eqref{E:integ1}, but
we will continue to refer to this extension to $\Rd$ by the expression of the integrand given above.
Making a change of variable $x=y^z+\delta v$, we get 
\begin{equation}\label{E:integ3}
h_\delta(z) = \int_{\Rd}\frac{f(y^z+\delta v)/G(y^z+\delta v) \ dv}
				{\left\{1+\delta^{-2}\drv{\gamma(y^z+\delta v)}{\gamma_j
									(y^z+\delta v)-z_j}\right\}^d}.
\end{equation} }
\medskip

Observe that 
\begin{equation}\label{E:junk}
\left|1+\delta^{-2}\drv{\gamma(y^z+\delta v)}{\gamma_j(y^z+\delta v)-z_j}\right| \ \geq \ 
	1+\delta^{-2}\er\left\{\drv{\gamma(y^z+\delta v)}{\gamma_j(y^z+\delta v)-z_j}\right\}.
\end{equation}
We compute
\begin{align}\label{E:compute}
\quad {} & \ \er\left\{\drv{\gamma(y^z+\delta v)}{\gamma_j(y^z+\delta v)-z_j}\right\} \\
       = & \ \er\left\{\drv{\gamma(y^z+\delta v)}{\gamma_j(y^z+\delta v)-\gamma(y^z)}+
			\drv{\gamma(y^z+\delta v)}{\gamma_j(y^z)-z_j}\right\} \notag \\
       = & \ \er\left\{\drv{\gamma(y^z+\delta v)}{\gamma_j(y^z+\delta v)-\gamma(y^z)}+
			\drv{\gamma(y^z)}{\gamma_j(y^z)-z_j}\right. \notag \\
\quad {} & \qquad \left. \add\left[\derrho{\gamma(y^z+\delta v)}-\derrho{\gamma(y^z)}\right]
								[\gamma_j(y^z)-z_j]\right\} \notag
\end{align}
From the fact that $\er\left\{\drv{\gamma(y^z)}{\gamma_j(y^z)-z_j}\right\}\geq 0 \ \forall z\in\overline\OM$,
and from Lemma \ref{L:pos} and \eqref{E:secOrdEst}, we have 
\begin{equation}\label{E:lobound1}
\er\left\{\drv{\gamma(y^z+\delta v)}{\gamma_j(y^z+\delta v)-z_j}\right\} \ 
		\geq \ \frac{C(\gamma)}{2}\ \delta^2|v|^2 \ \forall z\in\overline\OM\cap\mathcal{U}(\varrho).
\end{equation}
From \eqref{E:junk} and \eqref{E:lobound1} we have
\begin{equation}\label{E:unifbd*}
|h_\delta(z)| \ \leq \ \|f/G\|_\infty\int_{\Rd}\left\{1+\frac{C(\gamma)}{2}|v|^2\right\}^{-d} dv \quad
						\forall z\in\overline\OM\cap\mathcal{U}(\varrho).
\end{equation}
We note here that \cite[Lemma 2.4]{wR} makes it clear that $\|1/G\|_\infty<\infty$.
\smallskip

We now consider the case when $z\in\overline\OM\setminus\mathcal{U}(\varrho)$. 
Due to the fact that $\gamma[\BDC(0;\varrho)]\subsetneq\mathcal{U}(\varrho)$, we can find a uniform
constant $c^\prime>0$ such that
\begin{equation}\label{E:separation}
\er\left\{\drv{\gamma(x)}{\gamma_j(x)-z_j}\right\} \ \geq \ 
			c^\prime \ \er\left\{\drv{\gamma(x)}{\gamma_j(x)-\gamma_j(0)}\right\} \quad
				\forall x\in\BDO(0;\varrho).
\end{equation}
This time, we make a change of variable $x=\delta v$. For $z\in\overline\OM\setminus\mathcal{U}(\varrho)$,
this results in
\begin{equation}\label{E:integ3bis}
h_\delta(z) = \int_{\Rd}\frac{f(\delta v)/G(\delta v) \ dv}
				{\left\{1+\delta^{-2}\drv{\gamma(\delta v)}{\gamma_j
									(\delta v)-z_j}\right\}^d}.
\end{equation}
From Lemma \ref{L:pos} and \eqref{E:separation}
we can deduce that
\begin{equation}\label{E:unifbd**}
|h_\delta(z)| \ \leq \ \|f/G\|_\infty\int_{\Rd}\{1+c^\prime |v|^2\}^{-d} \ dv \quad
						\forall z\in\overline\OM\setminus\mathcal{U}{\varrho}.
\end{equation}
Claim (i) follows from \eqref{E:unifbd*} and \eqref{E:unifbd**}.
\medskip

The above argument actually yields the following observation, which we record.
\smallskip

\noindent{{\bf Fact. }\em There exists a uniform constant $\kappa>0$ such that the integrands occuring in 
\eqref{E:integ3} and \eqref{E:integ3bis} are dominated by the $\mathbb{L}^1$ function
\[
\Rd \ni v \mapsto \|f/G\|_\infty\{1+\kappa |v|^2\}^{-d}.
\] }
\medskip

\noindent{{\bf Claim (ii)} $\lim_{\delta\to 0} h_\delta(z)=0$ if
$z\in\overline\OM\setminus\gamma[\BDO(0;\varrho)]$. }
\smallskip

\noindent{Notice that if $z\in\overline\OM\setminus\gamma[\BDO(0;\varrho)]$, then
$\er\left\{\drv{\gamma(x)}{\gamma_j(x)-z_j}\right\} > 0 \ \forall x\in\BDO(0;\varrho)$. Thus
\[
\lim_{\delta\to 0}\frac{\delta^d f(x)/G(x)}
	{\left\{\delta^2+\drv{\gamma(x)}{\gamma_j(x)-z_j}\right\}^d}=0.
\]
In view of the fact recorded above, we can apply the dominated convergence theorem to \eqref{E:integ1}. This
results in Claim (ii).}
\medskip

\noindent{{\bf Claim (iii)} $\lim_{\delta\to 0} h_\delta[\gamma(s)]=f(s) \
\forall s\in\BDO(0;\varrho)$. }
\smallskip

\noindent{Refer to Lemma \ref{L:localmin}. When $z=\gamma(s)$ in that lemma, $y^z=\gamma(s)$. Equation
\eqref{E:integ3} reads as
\[
h_\delta[\gamma(s)] = \int_{\Rd}\frac{f(s+\delta v)/G(s+\delta v) \ dv}
				{\left\{1+\delta^{-2}\drv{\gamma(s+\delta v)}{\gamma_j
									(s+\delta v)-\gamma_j(s)}\right\}^d}.
\]
In view of \ref{E:keylimit}, the integrands occuring above tend to 
\[
\frac{f(s)}{G(s)}\left\{1+\frac{1}{2}\varQ \right\}^{-d}, \ 
\text{as $\delta\to 0$}.
\]
Once again, Claim (iii) follows from the dominated convergence theorem.}
\end{proof}
\bigskip

\section{The proof of Theorem \ref{T:main}}

For a $p\in\mn$, let $\strat$ stand for an arbitrary $d$-dimensional stratum, $d\geq 1$, of the 
local stratification \eqref{E:strata} of $\mn$ near $p$. Let $q\in\strat$ and let $\gamma : (\BDO(0;R),0)\to
(\strat,q)$ be a non-singular, real-analytic parametrization of $\strat$ near $q$. Notice that by the 
definition of $\strat$, $\gamma$ satisfies the hypotheses of Theorem \ref{T:big}. 
In view of Bishop's theorem (refer back to Theorem \ref{T:Bishop}), it would suffice to show that
for any compact $K\subset\BDO(0;\varrho)$ and any annihilating measure $\mu\perp A(\OM)$, $\mu[\gamma(K)]=0$, 
where $\varrho\equiv\varrho(\gamma)>0$ is the constant introduced in Theorem \ref{T:big}.
\medskip

Now, given a compact $K\subset\BDO(0;\varrho)$, let $\{\mathfrak{D}_\nu\}_{\nu\in\mathbb{N}}$ be a shrinking 
family of compact subsets such that
\begin{enumerate}
\item[(a)] $\mathfrak{D}_\nu\subset \BDO(0;\varrho)$,
\item[(b)] $\mathfrak{D}_{\nu+1}\subset \text{int}(\mathfrak{D}_\nu)$,
\item[(c)] $\cap_{\nu\in\mathbb{N}}=K$.
\end{enumerate}
Let $\chi_\nu\in\mathcal{C}^\infty(\BDO(0;\varrho);[0,1])$ be a bump function with
\[
\chi_\nu|_{\mathfrak{D}_{\nu+1}}\equiv 1, \quad\qquad \text{supp} \ \chi_\nu\subseteq\mathfrak{D}_\nu.
\]
\smallskip

We define $h^\nu_{\delta}\in A(\OM)$ by taking $f=\chi_\nu$ in the equation \eqref{E:integ1}. Let 
$\mu\perp A(\OM)$. 
By Theorem \ref{T:big} and the bounded convergence theorem, we have
\[
0 = \lim_{\delta\to 0}\int h^\nu_\delta d\mu = \int_{\gamma[\BDO(0;\varrho)]}\widetilde{\chi_\nu}d\mu,
\]
where $\widetilde{\chi_\nu}$ are given by the equations 
$\widetilde{\chi_\nu}[\gamma(x)]=\chi_\nu(x) \ \forall x\in\BDO(0;\varrho)$. Another passage to the limit
gives $\mu[\gamma(K)]=0$, and this is true for any $\mu\perp A(\OM)$. 
\medskip

We have just shown that each $\strat$ is a countable union of peak-interpolation sets for $A(\OM)$.
Since each of the finitely many points in $M_0$ ($M_0$ as given by \eqref{E:strata}) are peak points for
$A(\OM)$ (since $\OM$ is convex) , $\mn$ is a compact subset of $\bdy$ that is a countable union of
peak-interpolation sets for $A(\OM)$. Using Bishop's theorem again, we conclude that $\mn$ is a 
peak-interpolation set for $A(\OM)$.


\bigskip

\begin{thebibliography}{10}

\bibitem{eB}
E. Bishop, {\em A general Rudin-Carleson theorem}, Proc. Amer. Math. Soc. {\bf 13} (1962),
140-143.

\bibitem{BG}
T. Bloom and I. Graham, {\em On type conditions for generic real submanifolds in $\mathbb{C}^n$}, 
Invent. Math. {\bf 40} (1984) 217-43.

\bibitem{dC}
D. Catlin, {\em Boundary invariants of pseudoconvex domains}, Ann. of Math.(2) {\bf 120} (1984), 529-586.

\bibitem{jdA}
J.P. D'Angelo, {\em Real hypersurfaces, orders of contact, and applications}, Annals of Math. 
{\bf 115} (1982), 615-637.

\bibitem{HT}
G. M. Henkin and A.E. Tumanov, {\em Interpolation submanifolds of pseudoconvex manifolds}, Translations
Amer. Math. Soc. {\bf 115} (1980), 59-69.

\bibitem{sL}
S. {\L }ojasiewicz, {\em Sur le probl\`eme de la division}, Studia Math. {\bf 18} (1959), 87-136.

\bibitem{McN}
J.D. McNeal, {\em Convex domains of finite type}, J. Funct. Anal. {\bf 108} (1992), 361-373.

\bibitem{aN}
A. Nagel, {\em Smooth zero sets and interpolation sets for some algebras of holomorphic functions on
strictly pseudoconvex domains}, Duke Math. J. {\bf 43} (1976), 323-348.

\bibitem{NR}
A. Nagel and W. Rudin, {\em Local boundary behavior of bounded holomorphic functions}, Can. J. Math.
{\bf 30} (1978), 583-592.

\bibitem{wR}
W. Rudin, {\em Peak-interpolation sets of class $\mathcal{C}^1$}, Pacific J. Math. {\bf 75} (1978),
267-279.

\end{thebibliography}
\end{document}